\documentclass[reqno]{amsart} 
\usepackage[colorlinks]{hyperref}
\usepackage{amssymb}
\usepackage{latexsym}
\usepackage{amstext} 
\usepackage{color}
\usepackage{amsthm}  
\usepackage{amsmath}
\begin{document}
\newtheorem{theorem}{Theorem}    
\newtheorem{conjecture}{Conjecture}
\newtheorem{question}{Question}
\newtheorem{proposition}[theorem]{Proposition}
\newtheorem{lemma}[theorem]{Lemma}
\newtheorem{corollary}[theorem]{Corollary}
\newtheorem{example}[theorem]{Example}
\newtheorem{remark}[theorem]{Remark}
\newtheorem{definition}[theorem]{Definition}
\setlength{\parskip}{.15in} 
\setlength{\parindent}{0in} 

\title[\hfil  Conjectures in Number Theory]{Some conjectures in elementary  number theory\\ } 
\author[Angelo B. Mingarelli \hfilneg]{Angelo B. Mingarelli}   
\address{School of Mathematics and Statistics\\ 
Carleton University, Ottawa, Ontario, Canada, K1S\, 5B6}
\email[A. B. Mingarelli]{amingare@math.carleton.ca}

\date{February 21, 2013 }
\thanks{This research is partially supported by an NSERC Canada Discovery Grant}
\subjclass[2010]{}
\begin{abstract}
\noindent{We announce a number of conjectures associated with and arising from a study of primes and irrationals in $\mathbb{R}$. All are supported by numerical verification to the extent possible.}
\end{abstract}

\maketitle
\section*{The Conjectures}
\subsection*{Bhargava factorials}

For definitions and basic results dealing with Bhargava's factorial functions we refer to \cite{mb2}, \cite{mb3}, \cite{mb1} and \cite{chab}. Briefly, let $X \subseteq {\bf Z}$ be a finite or  infinite set of integers. Following \cite{mb1}, one can define the notion of a $p$-ordering on $X$ and use it to define a set of generalized factorials of the set $X$ inductively. By definition $0!_X=1$. Whenever $p$ a prime, we fix an element $a_0 \in X$ and, for $k \geq 1$, we select $a_k$ such that the highest power of $p$ dividing $\prod_{i=0}^{k-1}{(a_k-a_i)}$ is minimized. The resulting sequence of $a_i$ is then called a $p$-ordering of $X$. As one can gather from the definition, $p$-orderings are not unique, as one can vary $a_0$. On the other hand, associated with such a $p$-ordering of $X$ we define an associated $p$-sequence $\{\nu_k(X,p)\}_{k=1}^{\infty}$ by $$\nu_k(X,p) = w_p(\prod_{i=0}^{k-1}{(a_k-a_i)}),$$ where $w_p(a)$ is, by definition, the highest power of $p$ dividing $a$ (e.g., $w_2(80)=16$). One can show that although the $p$-ordering is not unique the associated $p$-sequence is independent of the $p$-ordering used. Since this quantity is an invariant it can be used to define generalized factorials of $X$ by setting \begin{equation}\label{mba} k!_X = \prod_{p}\nu_k(X,p),\end{equation} where the (necessarily finite) product extends over all primes $p$.

\begin{definition} \label{def2}\cite{af}. An abstract (or generalized) factorial is a function $!_{_{a}}: \mathbb{N} \to \mathbb{Z^+}$ that satisfies the following conditions:
\begin{enumerate}
\item $0!_{_{a}}=1$,
\item For every non-negative integers $n, k$, $0 \leq k \leq n$ the generalized binomial coefficients $$\binom{n}{k}_{_{a}} := \frac{n!_{_{a}}}{k!_{_{a}}(n-k)!_{_{a}}} \in \mathbb{Z^+},$$
\item For every positive integer $n$, $n!$ divides $n!_{_{a}}$.
\end{enumerate}
\end{definition}
It is easy to see that the collection of all abstract factorials forms a commutative semigroup under ordinary pointwise multiplication.  In fact, it is easy to see that Bhargava's factorial function is an abstract factorial. (Indeed, Hypothesis 1 of Definition~\ref{def2} is clear by definition of the factorial in question. Hypothesis 2 of Definition~\ref{def2} follows by the results in  \cite{mb1}. )

The context of these first three conjectures is the construction in \cite{mb1} as applied to the ring of integers. In this case, the factorial function for the set of rational primes 
$$\mathbb{P}=\{2,3,5,7,11, \dots\}$$ is given by \cite{mb1}
\begin{equation}\label{eqmbpr} n!_{_\mathbb{P}} = \prod_{p}^{} \displaystyle p^{\sum_{m=0}^{\infty}{[\frac{n-1}{p^m(p-1)}]}}.
\end{equation}
We call this simply the B-factorial for the set under consideration. In the sequel, the statement ``For every $n \geq 1$" means ``for every integer $n\geq 1$ for which the factorials are defined".

Let $\mathbb{P}_2 \subset \mathbb{P}$ denote the subset of all twin primes, i.e., those primes of the form $p, p+2$ as usual. Let $n!_{\mathbb{P}_2}$ denote the B-factorial of the set $\mathbb{P}_2$.
In the following conjectures the notation $w_p(n)$ is used to identify  the highest power of $p$ that divides $n$. So, for example, if $n$ has the representation  $n=2^{a_1}\alpha$ and $(\alpha,2)=1$, then $w_2(n) = 2^{a_1}$.
\begin{conjecture} For every $n \geq 1$,
$$\frac{n!_{\mathbb{P}_2}}{n!_{\mathbb{P}}} =  2\,w_2(n).$$

\end{conjecture}
In analogy with the preceding we let $\mathbb{P}_3 \subset \mathbb{P}$ denote that subset of all prime triplets of the form $p, p+2, p+6$.  Let $n!_{\mathbb{P}_3}$ denote the B-factorial of the set $\mathbb{P}_3$.
\begin{conjecture} For every $n \geq 1$,
$$\frac{n!_{\mathbb{P}_3}}{n!_{\mathbb{P}}} =  \left \{ \begin{array}{ll}
		3!\, w_2(n)w_3(n),	&	\mbox{if  \ \ $n$\ \ is\ \ even}, \\
		2 ,	&	\mbox{if  \ \ $n$\ \ is \ \ odd}. \\
		\end{array}
	\right. 
$$
\end{conjecture}

Next, let $\mathbb{P}_4 \subset \mathbb{P}$ denote that subset of all prime quadruplets written in the  form  $p, p+2, p+6, p+8$. Since $p, p+2$ and $p+6,p+8$ are both twin primes we can view $\mathbb{P}_4\subset \mathbb{P}_2$, and so we must have $n!_{\mathbb{P}_2}|n!_{\mathbb{P}_4}$, by [\cite{mb1}, Lemma~13]. In fact, we claim that,
\begin{conjecture} For  every $n \geq 1$,
$$\frac{n!_{\mathbb{P}_4}}{n!_{\mathbb{P}_2}} =  \left \{ \begin{array}{ll}
		3\,w_3(n),	&	\mbox{if  \ \ $n$\ \ is\ \ even}, \\
		1 ,	&	\mbox{if  \ \ $n$\ \ is \ \ odd}. \\
		\end{array}
	\right. 
$$
\end{conjecture}
These three conjectures have been verified using Crabbe's algorithm \cite{ac} to the limits available by the hardware. For motivation see \cite{af}.

\subsection*{Prime number inequalities}

Now let $p_n$ denote the n-th prime. Then, see \cite{af},
\begin{conjecture}
\begin{equation}\label{con}p_n \geq p_k + p_{n-k-1}, \quad\quad 1 \leq k \leq n-1,\end{equation} and all $n \geq 2$. 
\end{conjecture}
The validity of this conjecture implies that the function
$f : \mathbb{N} \to \mathbb{Z^+}$,
\[ f(n) = \left \{ \begin{array}{ll}
		1 ,	&	\mbox{if  \ \ $n = 0$}, \\
		1 ,	&	\mbox{if  \ \ $n =1 $}. \\
        p_{n-1}! ,	&	\mbox{if  \ \ $n \geq 2 $}. \\
		\end{array}
	\right. 
\]
 is an abstract factorial.  Thus, if true, it would follow from the results in \cite{af} that for any abstract factorial $n!_{_a}$, the quantity $\sum_{n\geq 1} 1/n!_{_a}f(n) \notin \mathbb{Q}$.

\subsection*{Ap\'ery numbers}

We define the Ap\'ery numbers $A_n, B_n$ recursively, as usual, by setting $A_0=1, A_1=5$; $B_0=0, B_1=6$ whose  general tems are given by the recurrence relations
$$A_{n+1} = (P(n)A_n-n^3A_{n-1})/(n+1)^3,$$
and
$$B_{n+1} = (P(n)B_n-n^3B_{n-1})/(n+1)^3,$$
where $P(n)$ is the polynomial 
$$P(n) = 34n^3+51n^2+27n+5.$$
In a singular argument Ap\'ery \cite{rap} showed that $B_n/A_n \to \zeta(3)$ as $n \to \infty$ where $\zeta$ is the usual Riemann zeta function. In addition, he proved that $\zeta(3)$ is irrational (though no explicit formula akin to the one known for the values of $\zeta$ at positive even integers was given). More explicit proofs appeared since, e.g., \cite{fb}, \cite{vdp}, \cite{hc1} among others. (See \cite{am} for extensions of the series acceleration method found in  {\rm [Fischler~\cite{sf}, Remarque 1.3]} to integer powers of $\zeta(3)$.)

Here we propose using an old irrationality criterion due to Brun \cite{vb} (see also \cite{vb2}) in order to formulate a conjecture that, if true, would give another proof of the irrationality of $\zeta(3)$. Let $x_n$ be a sequence of real numbers and $\Delta$ the forward difference operator defined by $\Delta x_n = x_{n+1}-x_n$. 

\begin{theorem}\label{brun}(Brun, \cite{vb}) Let $x_n \in \mathbb{Z^+}$ be an increasing sequence and $y_n \in \mathbb{Z^+}$ be such that $\Delta (y_n/x_n) > 0$. If
\begin{equation}\label{del}
\delta_n \equiv \Delta\left ( \Delta y_n/\Delta x_n \right) < 0,
\end{equation}
then $y_n/x_n$ converges to an irrational number.
\end{theorem}
Although Brun claimed later \cite{vb2} that ``\ldots this theorem is simple but unfortunately not very useful" we show that perhaps it may be used to prove the irrationality of $\zeta(3)$. 

The idea is as follows: It is known that the sequence $A_n$ of Ap\'ery numbers is an increasing sequence of positive integers \cite{hc1} and although the $B_n$ is not necessarily a sequence of integers, the weighted sequence $e_nB_n$ is such a sequence where $e_n = 2\cdot ({\rm lcm}\{1,2,\ldots,n\})^3$, \cite{hc1}. In addition, the sequence $B_n/A_n=e_nB_n/e_nA_n$ is increasing, \cite{hc1} and it is easily proved that the sequence $e_nA_n$ is increasing as well.

Thus, setting $x_n=e_nA_n$ and $y_n=e_nB_n$ we see that the requirements $x_n$ is increasing and $y_n/x_n$ increasing are met in Theorem~\ref{brun} (all sequences being positive and all integers). We anticipate the following
\begin{conjecture} There is an unbounded subsequence of positive integers $n_k \to \infty$ such that 
$\delta_{n_k}  < 0.$
\end{conjecture}
Since it is known that $y_{n}/x_n$ increases to $\zeta(3)$, clearly $y_{n_k}/x_{n_k}$ does the same for any subsequence. Hence, an affirmative answer to the previous conjecture implies the irrationality of $\zeta(3)$ by Brun's irrationality theorem,Theorem~\ref{brun}. The numerical evidence seems to point to a stronger conjecture however. Indeed, it appears as if
\begin{conjecture} For every integer $N\geq 2$, there is an $n \in \mathbb{Z^+}$ such that all
$$\delta_n, \delta_{n+1}, \delta_{n+2}, \ldots, \delta_{n+N} < 0.$$
\end{conjecture}
Of course this result, if true, implies the previous conjecture.


\begin{thebibliography}{+99}
\bibitem {rap} R. Ap\'{e}ry, {\it Irrationalit\'{e} de $\zeta{2}$ et $\zeta{3}$,} Ast\'{e}risque, {\bf 61} (1979), 11-13
\bibitem{fb} F. Beukers, {\it A note on the irrationality of $\zeta(2)$ and $\zeta(3)$}, Bull. London Math. Soc., {\bf 11} (1979), 268-272
\bibitem{mb2} M.Bhargava, {\it P-orderings and polynomial functions on arbitrary subsets of Dedekind rings}, J. reine angew. Math., {\bf 490} (1997), 101-127
\bibitem{mb3} M.Bhargava, {\it Generalized factorials and fixed divisors over subsets of Dedekind domains}, J. Number Theory, {\bf 72} (1) (1998), 67-75
\bibitem{mb1} M.Bhargava, {\it The factorial function and generalizations}, Amer. Math. Monthly, {\bf 107} (2000), 783-799
\bibitem{vb} V. Brun, {\it Ein Satz \"uber Irrationalit\"at}, Arch. for Math. og Naturvidenskab
(Kristiania) {\bf 31} (1910), 3-6.
\bibitem{vb2} V. Brun and F. F. Knudsen, {\it On the possibility of finding certain criteria for the irrationality of a number defined as a limit of a sequence of rational numbers}, Math. Scand., {\bf 31} (1972), 231-236.
\bibitem{chab} J.-L. Chabert and P.-J. Cahen, {\it Old problems and new questions around integer-valued polynomials and factorial sequences}, in {\em Multiplicative Ideal Theory in Commutative Algebra}, J. Brewer {\it et al} eds., Springer, New York (2006),  89-108.
\bibitem{hc1} H. Cohen, {\it D\'{e}monstration de l'irrationalit\'{e} de $\zeta(3)$}, S\'{e}minaire de Th\'{e}orie des Nombres, Grenoble, 5 octobre 1978, No. 6, 9p.
\bibitem{ac} A. M. Crabbe, {\it Generalized Factorial Functions and Binomial Coefficients}, Undergraduate Honors Thesis, Trinity University, USA, (2001), 35 pp.
\bibitem{sf} S. Fischler, {\it Irrationalit\'{e} de valeurs de z\^{e}ta}, S\'{e}minaire Bourbaki, 55\`{e}me ann\'{e}e, 2002-2003, No. 910, Nov. 2002.
\bibitem{af} A. B. Mingarelli, {\it Abstract Factorials}, submitted.
\bibitem{am} A. B. Mingarelli, {\it On a discrete version of a theorem of Clausen and its applications}, To appear in Acta Math.  Acad.  Paedagogicae  Nyíregyháziensis, {\bf 29} (1), (2013).
\bibitem{vdp} A. Van der Poorten, {\it A proof that Euler missed: Ap\'{e}ry's proof of the irrationality of $\zeta (3)$}, Math. Intelligencer, {\bf 1} (4) (1979), 195-203
\end{thebibliography}
\end{document}